\documentclass[12pt]{amsart}
\usepackage{amsfonts, amsmath, latexsym, vmargin, epsfig}
\usepackage{epsf}
\usepackage{url}
\usepackage[ps2pdf, colorlinks=true,urlcolor=blue]{hyperref}

\title{Hypercube embedding of Wythoffians}

\author{Michel Deza}
\address{Michel Deza, LIGA, \'Ecole Normale Sup\'erieure, 45 rue d'Ulm, 75005 Paris, France}
\email{Michel.Deza@ens.fr}

\author{Mathieu Dutour Sikiri\'c}
\address{Mathieu Dutour Sikiri\'c, Rudjer Bo\u skovi\'c Institute, Bijenicka 54, 10000 Zagreb, Croatia}
\email{mdsikir@irb.hr}

\author{Sergey Shpectorov}
\address{School of Mathematics, University of Birmingham, Edgbaston,
Birmingham, B15 2TT, United Kingdom}
\email{sergeys@for.mat.bham.ac.uk}
\thanks{The first author was partly supported by JAIST, Japan.
The second author was partly supported by the Croatian Ministry of Science, Education and Sport under contract 098-0982705-2707.
The third author was partly supported by an NSA grant.}

\newtheorem{proposition}{Proposition}

\newtheorem{conjecture}{Conjecture}

\newcommand{\RR}{\ensuremath{\mathbb{R}}}

\newcommand{\ZZ}{\ensuremath{\mathbb{Z}}}

\newcommand{\cK}{\ensuremath{{\mathcal K}}}
\newcommand{\cM}{\ensuremath{{\mathcal M}}}

\newcommand{\Dl}{\Delta}
\newcommand{\Gm}{\Gamma}
\newcommand{\Om}{\Omega}

\newcommand{\gm}{\ensuremath{\gamma}}

\newcommand{\eps}{\ensuremath{\varepsilon}}
\newcommand{\lm}{\ensuremath{\lambda}}

\newcommand{\Med}{\ensuremath{\rm Med}}

\begin{document}
\maketitle

\begin{abstract}
The Wythoff construction takes a $d$-dimensional polytope $P$, a subset $S$
of $\{0,\dots, d\}$ and returns another $d$-dimensional polytope $P(S)$.
If $P$ is a regular polytope, then $P(S)$ is vertex-transitive.
This construction builds a large part of the Archimedean polytopes
and tilings in dimension $3$ and $4$.

We want to determine, which of those Wythoffians $P(S)$ with regular $P$
have their skeleton or dual skeleton isometrically embeddable
into the hypercubes $H_m$ and half-cubes $\frac{1}{2}H_m$.
We find six infinite series, which, we conjecture, cover all cases for
dimension $d>5$ and some sporadic cases in dimension $3$ and $4$
(see Tables \ref{WythoffEmbeddable3} and \ref{WythoffEmbeddable4}).

Three out of those six infinite series are explained by a general
result about the embedding of Wythoff construction for Coxeter groups.
In the last section, we consider the Euclidean case; also,
zonotopality of embeddable $P(S)$ are addressed throughout the text.

\end{abstract}

\section{Wythoff kaleidoscope construction}\label{WythoffSection}

A {\em flag} in a poset is an arbitrary completely ordered subset. We say
that a connected poset $\cK$ is a {\em $d$-dimensional complex} (or,
simply, a {\em $d$-complex}) if every maximal flag in $\cK$ has size
$d+1$. In a $d$-complex $\cK$ every element $x$ can be uniquely assigned a
number $\dim(x)\in\{0,\ldots,d\}$, called the {\em dimension} of $x$, in
such a way, that the minimal elements of $\cK$ have dimension zero and
$\dim(y)=\dim(x)+1$ whenever $x<y$ and there is no $z$ with $x<z<y$.

The elements of a complex $\cK$ are called {\em faces}, or {\em $k$-faces}
if the dimension of the face needs to be specified.  Furthermore, $0$-faces
are called {\em vertices} and $d$-faces (maximal faces) are called {\em
facets}. If we reverse the order on $\cK$ then the resulting poset
$\cK^\ast$ is again a $d$-complex, called the {\em dual complex}. Clearly,
the vertices of $\cK^\ast$ are the facets of $\cK$ and, more generally, the
dimension function on $\cK^\ast$ is given by $\dim^\ast(x)=d-\dim(x)$.

We will often use the customary geometric language. If $x<y$ and
$\dim(x)=k$, we will say that {\em $x$ is a $k$-face of $y$}.

A $d$-complex is a {\em polytope} if every submaximal flag (that is, a flag
of size $d$) is contained in exactly two maximal flags.  In the polytopal
case, $1$-faces are called {\em edges}, because each of them has exactly two
vertices. Starting from the next section we will deal exclusively with
polytopes. The {\em skeleton} of a polytope $\cK$ is the graph formed by
all vertices and edges of $\cK$.

For a flag $F\subset\cK$ define its {\em type} as the set
$t(F)=\{\dim(x)|x\in F\}$. Clearly, $t(F)$ is a subset of
$\Dl=\{0,\ldots,d\}$ and, reversely, every subset of $\Dl$ is the type of
some flag.

Let $\Om$ be the set of all nonempty subsets of $\Dl$ and fix an
arbitrary $V\in\Om$. For two subsets $U,U'\in\Om$ we say that $U'$
{\em blocks} $U$ (from $V$) if for all $u\in U$ and $v\in V$ there is
a $u'\in U'$, such that $u\le u'\le v$ or $u\ge u'\ge v$. This defines
a binary relation on $\Om$, which we will denote as $U'\le U$. We also
write $U'\sim U$ if $U'\le U$ and $U\le U'$, and we write $U'<U$ if
$U'\le U$ and $U\not\le U'$.

It is easy to see that $\le$ is reflexive and transitive, which implies
that $\sim$ is an equivalence relation. Let $[U]$ denote the equivalence
class containing $U$. It will be convenient for us to choose canonic
representatives in equivalence classes. It can be shown that if $U\sim U'$
then $U\cap U'\sim U\sim U\cup U'$. This yields that every equivalence
class $X$ contains a unique smallest (under inclusion) subset $m(X)$ and
unique largest subset $M(X)$. If $X=[U]$ then $m(X)$ and $M(X)$ can be
specified as follows: $m(X)$ is the smallest subset of $U$ that blocks $U$,
while $M(X)$ is the largest subset of $\Dl$ that is blocked by $U$. The
subsets $m(X)$ will be called the {\em essential} subsets of $\Dl$ (with
respect to $V$). Let $E=E(V)$ be the set of all essential subsets of
$\Dl$. Clearly, the above relation $<$ is a partial order on $E$. Also,
$V\in E$ and $V$ is the smallest element of $E$ with respect to $<$.

We are now ready to explain the {\em Wythoff construction}. Naturally, our
description is equivalent to the one given in \cite{Cox35} and
\cite{Cox73}, that generalized the original paper \cite{Wyt}.
Other relevant reference to the subject are \cite{mcmullen} and
\cite{Sch90}.
Suppose $\cK$ is a $d$-complex and let $\Dl$,
$\Om$, $V$, $\le$, and $E$ be as above. The {\em Wythoff complex} (or {\em
Wythoffian}) $\cK(V)$ consists of all flags $F$ such that $t(F)\in E$. For
two such flags $F$ and $F'$, we have $F'<F$ whenever $t(F')<t(F)$ and $F'$
is {\em compatible} with $F$ (that is, $F\cup F'$ is a flag). It can be
shown that $\cK(V)$ is again a $d$-complex and that
$\dim(F)=d+1-|M([t(F)])|$. It can also be shown that if $\cK$ is a
$d$-polytope then $\cK(V)$ is again a $d$-polytope for all $V$.
For a concrete Euclidean realization of such polytopes,
see \cite{HE93}.

Since there are $2^{d+1}-1$ different subsets $V$, there are, in general,
$2^{d+1}-1$ different Wythoffians constructed from the same complex
$\cK$. It is easy to see that $\cK(V)=\cK^\ast(d-V)$, where $d-V=\{d-v|v\in
V\}$. This means that the dual complex does not produce new
Wythoffians. Furthermore, in the case of self-dual complexes (that is, where
$\cK\cong\cK^\ast$), this reduces the number of potentially pairwise
non-isomorphic Wythoffians to $2^{d}+2^{\lceil
\frac{d-1}{2}\rceil}-1$.

Some of the Wythoffians are, in fact, familiar complexes. First of all,
$\cK(\{0\})=\cK$ and $\cK(\{d\})=\cK^\ast$. Furthermore, $\cK(\{1\})$ is
also known as the {\em median} complex $\Med(\cK)$ of $\cK$ and the dual of
$\cK(\Dl)$ is known as the {\em order complex} of $\cK$ (see
\cite{RS}). We will call $\cK(\Dl)$ the {\em flag complex} of $\cK$.
Thus, the order complex is the dual of the flag complex.

Since in this paper we are going to deal with the skeletons of $\cK(V)$ and
$\cK(V)^\ast$ (in the polytopal case), we need to understand elements of
$\cK(V)$ of types $0$, $1$, $d-1$, and $d$. Since $V$ is the unique
smallest essential subset, the vertices ($0$-faces) of $\cK(V)$ are the
flags of type $V$. For a flag $F$ to be a $1$-face of $\cK(V)$, $U=t(F)$
must have the property that $M([U])$ misses just one dimension $k$ from
$\Dl$. Clearly, $k$ must be in $V$. Now $U=U_k$ can be readily
computed. Namely, $U_k$ is obtained from $V$ by removing $k$ and including
instead the neighbors of $k$ (that is, $k-1$ and/or $k+1$).  Thus, $\cK(V)$
has exactly $|V|$ types of $1$-faces. Turning to the facets ($d$-faces), we
see that, for $F$ to be a facet of $\cK(V)$, we need that $U=t(F)$ be an
essential subset of size one, such that $M([U])=U$. The latter condition
can be restated as follows: $U$ should block no other $1$-element set. From
this we easily obtain that the relevant sets $U=\{k\}$ are those for which
$k=0$ (unless $V=\{0\}$), $k=d$ (unless $V=\{d\}$), or
$\min(V)<k<\max(V)$. Finally, if $F$ is a $(d-1)$-face then $U=t(F)$ is
essential of size one or two and $M([U])$ is of size exactly two. We will
not try to make here a general statement about all such subsets
$U$. However, in the concrete situations below, it will be easy to list
them all.

\section{Archimedean Wythoffians: $d=3$, $4$} \label{dimension 2 and 3}

In this section we start looking at particular examples of Wythoffians,
namely, at the {\em Archimedean Wythoffians}. These polytopes come by the
Wythoff construction from the regular convex polytopes. A complex (in
particular, a polytope) is called {\em regular} if its group of symmetries
acts transitively on the set of maximal flags. {\em Convex polytopes} are
the ones derived from convex hulls $H$ of finite sets of points in
$\RR^d$. (We assume that the initial set of points contains $d+1$ points in
general position; equivalently, the interior of $H$ is nonempty.) The faces
of the polytope are the convex intersections of the boundary of $H$ with
proper affine subspaces of $\RR^d$. In particular, the polytope is
$(d-1)$-dimensional, rather than $d$-dimensional.

It is well-known that the regular convex polytopes
fall into four infinite series: regular $p$-gon,
simplices $\alpha_d$, hyperoctahedra $\beta_d$,
and hypercubes $\gm_d$; and five sporadic examples: the icosahedron and
dodecahedron for $d=3$, and the $24$-cell, $600$-cell, and $120$-cell for
$d=4$.

The {\em half-cube graph} $\frac{1}{2}H_m$
(respectively {\em Johnson graph} $J(m,n)$) is the graph formed by all
$x\in \{0,1\}^m$ such that $\sum_{i=1}^m x_i$ is even
(respectively equal $n$) with two vertices
adjacent if their Hamming distance is $2$.

We are interested in the following 

\medskip\noindent
{\bf Main Question:} {\em Which Archimedean Wythoffians have skeleton graph
or dual skeleton graph isometrically embeddable, for a suitable $m$, in the
hypercube graph $H_m$ or half-cube graph $\frac{1}{2}H_m$?} 

\medskip
Recall (\cite{AD80,GW,DS,Shp,DS1} and books \cite{book1,book2})
that a mapping $\phi$ from a graph $\Gm$ to a graph $\Gm'$ is an
{\em isometric embedding} if $d_{\Gm'}(\phi(u),\phi(v))=d_\Gm(u,v)$ for all
$u,v\in\Gm$. For brevity, we will often shorten ``isometric embedding'' to
just ``embedding''. Notice that $H_m$ is an isometric subgraph of
$\frac{1}{2}H_{2m}$, which means that every graph isometrically embeddable
in a hypercube is also embeddable in a half-cube.
There is also an intermediate class of graphs---those that are
embeddable in a Johnson graph $J(m,n)$.
A graph is said to be {\em hypermetric} if its path-metric satisfies
the inequality
\begin{equation*}
\sum_{1\leq i<j\leq n} b_ib_j d_{G}(i,j)\leq 0
\end{equation*}
for any vector $b\in \ZZ^n$ with $\sum_{i}b_i=1$.
In the special case, when $b$ is a permutation of
$(1,1,1,-1,-1,0,\dots,0)$, the above inequality is called {\em $5$-gonal}.
The validity of hypermetric inequalities is necessary for embeddability
but not sufficient: an example of hypermetric, but not embeddable
graph $K_7-C_5$ (amongst those, given in Chapter 17 of \cite{book1}).

Below, when we state our results on embeddability of the skeleton
graphs $\Gm$, we will indicate the smallest class in the above hierarchy,
containing $\Gm$.

We remark that $\gm_d$ is dual to $\beta_d$, which means that they
produce the same Wythoffians. Thus, one can skip the case $\cK=\gm_d$
altogether.  Similarly, we can skip the cases where $\cK$ is the
dodecahedron, since the latter is dual to the icosahedron, and the
$120$-cell, since it is dual to the $600$-cell.

In the remainder of this section we state the results of a computer
calculation carried out in the computer algebra system {\tt GAP}
\cite{GAP}.

We start with the case $d=3$. In this case $\cK$ is $2$-dimensional,
that is, $\cK$ is a {\em map} (and so, we switch to the notation
$\cM=\cK$). It is easy to see that $\cM(V)$ with $V$=$\{0\}$,
$\{0,1\}$, $\{0,1,2\}$, $\{0,2\}$, $\{1\}$, $\{1,2\}$, and $\{2\}$
correspond, respectively, to the following maps: original map $\cM$,
truncated $\cM$, truncated $\Med(\cM)$, $\Med(\Med(\cM))$, 
$\Med(\cM)$, truncated $\cM^\ast$ and $\cM^\ast$.

In Table \ref{WythoffEmbeddable3} we give a complete answer to our Main
Question in the case $d=3$. The table lists all Archimedean Wythoffians and
dual Wythoffians, whose skeleton graph is embeddable. The details of the
embedding, such as the dimension of the embedding and whether or not it is
equicut, are also provided. Recall that an embedding of a graph $\Gm$ in a
hypercube is called {\em equicut} if each cut on $\Gm$, produced by a
coordinate of the hypercube, splits $\Gm$ in half. An embedding is called
{\em $q$-balanced} if each coordinate cut on $\Gm$ has parts of sizes $q$
and $|\Gm|-q$. We will indicate in the table whether the embedding is
equicut, $q$-balanced, or neither. Finally, for brevity, we truncated in
the table the word ``truncated'' to just ``tr''.
\begin{table}[ht]
\begin{center}
\begin{tabular}{||c|c|c|c||}
\hline
\hline
Embeddable Wythoffian                      & $n$ & embedding & equicut?\\
\hline
Tetrahedron$=\alpha_3(\{0\})=\alpha_3(\{2\})$    & 4 & $=J(4,1)$; $=\frac{1}{2}H_3$ & $q=1$; yes\\
Octahedron$=\beta_3(\{0\})=\alpha_3(\{1\})$     & 6 & $=J(4,2)$ & yes\\
Cube$=\beta_3(\{2\})=\beta_3(\{0\})^\ast$      & 8 & $=H_3$ & yes\\
Icosahedron$=Ico(\{0\})$                   & 12 & $\frac{1}{2}H_6$ & yes\\
Dodecahedron$=Ico(\{2\})$                  & 20 & $\frac{1}{2}H_{10}$ & yes\\
\hline
(tr Tetrahedron)$^\ast=\alpha_3(\{0,1\})^\ast=\alpha_3(\{1,2\})^\ast$ & 8 & $\frac{1}{2}H_7$ & no\\
(Cuboctahedron)$^\ast=\beta_3(\{1\})^\ast=\alpha_3(\{0,2\})^\ast$ & 14 & $H_4$ & yes\\
(tr Cube)$^\ast=\beta_3(\{1,2\})^\ast$       & 14 & $J(12,6)$ & no\\
Rhombicuboctahedron$=\beta_3(\{0,2\})$       & 24 & $J(10,5)$ & yes\\
tr Octahedron$=\beta_3(\{0,1\})=\alpha_3(\{0,1,2\})$ & 24 & $H_6$ & yes\\
(tr Icosahedron)$^\ast=Ico(\{0,1\})^\ast$  & 32 & $\frac{1}{2}H_{10}$ & yes\\
(Icosidodecahedron)$^\ast=Ico(\{1\})^\ast$ & 32 & $H_6$ & yes\\
(tr Dodecahedron)$^\ast=Ico(\{1,2\})^\ast$ & 32 & $\frac{1}{2}H_{26}$ & no\\
tr Cuboctahedron$=\beta_3(\{0,1,2\})$        & 48 & $H_9$ & yes\\
Rhombicosidodecahedron$=Ico(\{0,2\})$      & 60 & $\frac{1}{2}H_{16}$ & yes\\
tr Icosidodecahedron$=Ico(\{0,1,2\})$      & 120 & $H_{15}$ & yes\\
\hline
\hline
\end{tabular}
\end{center}
\caption{Embeddable Archimedean Wythoffians for $d=3$.}
\label{WythoffEmbeddable3}
\end{table}

A striking property of this table is that it contains all possible
Wythoffians (all five regular polytopes and 11 of the 13 Archimedean
polytopes; missing are the Snub Cube and Snub Dodecahedron, which are not
Wythoffian). Furthermore, for each of these polytopes, exactly one of the
skeleton and the dual skeleton is embeddable.

This nice picture does not extend to the case $d=4$, where far fewer
embeddings exist. Our Table \ref{WythoffEmbeddable4} gives a
complete answer to the Main Question. 
\begin{table}[ht]
\begin{center}
\begin{tabular}{||c|c|c|c||}
\hline
\hline
Embeddable Wythoffian             & $n$ & embedding & equicut?\\
\hline
$\alpha_4=\alpha_4(\{0\})=\alpha_4(\{3\})$ & 5 & $=J(5,1)$ & $q=1$\\
$\beta_4=\beta_4(\{0\})$              & 8 & $=\frac{1}{2}H_{4}$ & yes\\
$\gm_4=\beta_4(\{3\})=\beta_4(\{0\})^\ast$ & 16 &  $=H_4$ & yes\\
\hline
$\alpha_4(\{1\})=\alpha_4(\{2\})$       & 10 & $=J(5,2)$ & $q=4$\\
$\alpha_4(\{0,3\})^\ast$             & 30 & $H_5$ & yes\\
$\beta_4(\{0,3\})$                  & 64 & $\frac{1}{2}H_{12}$ & yes\\
$\alpha_4(\{0,1,2,3\})$              & 120 & $H_{10}$ & yes\\
$\beta_4(\{0,1,2\})=24-cell(\{0,1\})=24-cell(\{2,3\})$ & 192 & $H_{12}$ & yes\\
$\beta_4(\{0,1,2,3\})$              & 384 & $H_{16}$ & yes\\
$24-cell(\{0,1,2,3\})$            & 1152 & $H_{24}$ & yes\\
$600-cell(\{0,1,2,3\})$           & 14400 & $H_{60}$ & yes\\
\hline
\hline
\end{tabular}
\end{center}
\caption{Embeddable Archimedean Wythoffians for $d=4$.}
\label{WythoffEmbeddable4}
\end{table}

Notice that the total number of Archimedean Wythoffians for $d=4$ is
$45$, see \cite{Conway,moller1,moller2}.
Thus, Table \ref{WythoffEmbeddable4} indicates that the
embeddable cases become more rare as $d$ grows, and that, likely, there are
only finitely many infinite series of embeddings.  Furthermore, Tables
\ref{WythoffEmbeddable3} and \ref{WythoffEmbeddable4} lead us to a number 
of concrete conjectures about possible infinite series of embeddings. In
the next section we resolve those conjectures in affirmative by
constructing the series and verifying the embedding properties.

We conclude this section with some further remarks about the embeddings in
Tables \ref{WythoffEmbeddable3} and \ref{WythoffEmbeddable4}. The majority
of these embeddings are unique. The only exception is the Tetrahedron
$\alpha_3$, whose skeleton, the complete graph $K_4$, has two isometric
embeddings.  We also checked that all the skeleton graphs for $d=3$, that
turn out to be non-embeddable, violate, moreover, the so-called {\em
$5$-gonal inequality}.

\section{Relation with Coxeter groups}\label{SectionCoxeter}

A group $\mathsf{W}$ is a {\em Coxeter group} if it is a group generated by
a set $S=\{s_0, \dots, s_{d-1}\}$ with the elements $s_i$ satisfying to
the relations $s_i^2=1$ and $(s_is_j)^{m_{ij}}=1$ with $m_{ij} \geq 2$.
We refer to \cite{humphreyscoxeter} for all facts used in this Section.

One can encode the matrix $\mathsf{W}$ by a {\em Coxeter graph} on vertices
$\{0, \dots, d-1\}$ with two edges being adjacent if $m_{ij}\geq 3$.

The group is called {\em irreducible} if the Coxeter graph is
connected. The irreducible finite Coxeter groups are classified and
denoted by $\mathsf{A}_d$, $\mathsf{B}_d$, $\mathsf{D}_d$, $\mathsf{E}_6$, $\mathsf{E}_7$, $\mathsf{E}_8$, $\mathsf{F}_4$, $\mathsf{H}_3$,
$\mathsf{H}_4$, $\mathsf{I}_2(p)$ (see \cite[page 34]{humphreyscoxeter}).
A Coxeter group is the symmetry group of a regular polytope
if and only if its Coxeter graph is a path.
A finite Coxeter group $\mathsf{W}$ can be represented as a group of isometries
of $\RR^d$ with the $s_i$ being reflexions and their hyperplanes
$H_i$ delimiting a fundamental domain ${\mathcal S}$.

We now define the Wythoff construction for Coxeter groups.
An algebraic way is explained in \cite{maxwell} but we choose
an easier and more geometric way following \cite{patera}.
Take a finite Coxeter group $\mathsf{W}$ with a fundamental simplex ${\mathcal S}$.
If $T\subseteq \Delta=\{0,\dots, d-1\}$, then we take a point $v\in {\mathcal S}$
such that $v\in H_i$ if and only if $i\notin T$.
The Wythoff construction $\mathsf{W}(T)$ is then defined as the convex hull
of the orbit $\mathsf{W}(v)$. The orbit $\mathsf{W}(v)$ is on the sphere and in
\cite{patera} it is proved, using Delaunay polytopes, that the
combinatorics of $\mathsf{W}(T)$ depends only on $\mathsf{W}$ and $T$.

If $T=\{0, \dots, d-1\}$, then $\mathsf{W}(T)$ is the flag
complex of $\mathsf{W}$.
Its skeleton is the {\em Cayley graph} $Cay(\mathsf{W}, S)$
of $\mathsf{W}$ obtained from the canonical generating
set $S=\{s_0, \dots, s_{d-1}\}$.

Given a Coxeter group $\mathsf{W}$, denote by $T$ the set of elements, which
are conjugate to an element of $S$.

\begin{proposition}\label{ProbablyFolflore}
If $\mathsf{W}$ is a finite Coxeter group, $S$ is its canonical
generating set,
then the Cayley graph $Cay(\mathsf{W}, S)$ is isometrically
embeddable into $H_{|T|}$.
\end{proposition}
\proof Any finite Coxeter group can be realized as a group of isometries
of a space $\RR^d$. The elements of $T$ are realized as reflections
and the corresponding hyperplanes realize a plane arrangement and $\mathsf{W}$
acts simply transitively on its cells.
By a theorem of Eppstein (\cite{epstein}), the graph formed by the
cells with two cells adjacent if they share a $(d-1)$-dimensional
face, is isometrically embeddable into the hypercube $\{0,1\}^M$
with $M$ being the number of hyperplanes, i.e. $|T|$. \qed

See in Table \ref{RegularPolytopesAndCoxeterGroups} the dimensions
of the hypercube of embeddings with the
corresponding regular polytopes if existing.

\begin{table}\label{RegularPolytopesAndCoxeterGroups}
\begin{center}
\begin{tabular}{||c|c|c||}
\hline
$\mathsf{W}$      &  $|T|$     & regular polytope\\
\hline
$\mathsf{A}_d$    &  $d(d+1)/2$  &  $\alpha_d$\\
$\mathsf{B}_d$    &  $d^2$       &  $\gamma_d$ or $\beta_d$\\
$\mathsf{D}_d$    &  $d(d-1)$    &  none\\
$\mathsf{E}_6$    &  $36$        &  none\\
$\mathsf{E}_7$    &  $63$        &  none\\
$\mathsf{E}_8$    &  $120$       &  none\\
$\mathsf{F}_4$    &  $24$        &  $24$-cell\\
$\mathsf{H}_3$    &  $15$        &  Dodecahedron or Icosahedron\\
$\mathsf{H}_4$    &  $60$        &  $120$-cell or $600$-cell\\
$\mathsf{I}_2(p)$ &  $p$         &  regular $p$-gon\\
\hline
\end{tabular}
\end{center}
\caption{Embeddings of flag complexes of finite irreducible Coxeter groups.}
\end{table}

This proposition is, certainly, folklore, but we could not
find an appropriate reference. For example, \cite{mirror}
address the $3$-dimensional case.
Also, in \cite{niblo}, an embedding of finitely generated Coxeter group
into cube complexes is given. But this is a topological embedding, while
our embeddings are isometric.

The above setting can be generalized to affine Coxeter groups $\widetilde{\mathsf{W}}$.
An {\em affine Coxeter group} is a group $\widetilde{\mathsf{W}}$ obtained by adding a
reflection $s_d$ along an hyperplane not passing though $0$ to a
finite Coxeter group $\mathsf{W}$.
The group $\widetilde{\mathsf{W}}$ then has a fundamental domain delimited
by $d+1$ hyperplanes $H_0, \dots, H_{d}$.
If $\widetilde{T}$ denotes the set of classes of
parallel hyperplanes of $\widetilde{\mathsf{W}}$,
then its Cayley graph $Cay(\mathsf{W}, S)$ embeds
into $\ZZ^{|\widetilde{T}|}$ and one has $|\widetilde{T}|=|T|$
according to the nomenclature in \cite{humphreyscoxeter}.

Note also that (see \cite{zieglerreiner,Cammond}),
given a finite Coxeter group of root system $R=\{r_1, \dots, r_N\}$,
the zonotopal polytope $[-r_1, r_1]+\dots + [-r_N, r_N]$ is, actually,
the Wythoff construction of $\mathsf{W}$ on $S$.
It is easy to see that this zonotopal embedding is isometric
embedding into $H_{|T|}$ given in Proposition \ref{ProbablyFolflore}.

\section{Infinite series of embeddings}

Since in this section we are only interested in the infinite series of
embeddings, we restrict ourselves to the cases $\cK=\alpha_d$ and
$\beta_d$. These polytopes can be described in combinatorial terms as
follows: The faces of $\alpha_d$ are all proper nonempty subsets of the
set $\{1,\ldots,d+1\}$. The order on $\alpha_d$ is given by containment,
and the dimension of the face $X$ is $|X|-1$. Clearly, $\alpha_d$ has
${d+1\choose k+1}$ faces of dimension $k$. The faces of $\beta_d$ are
the sets $\{\pm i_1,\ldots,\pm i_k\}$, where the signs are arbitrary
and $\{i_1,\ldots,i_k\}$ is a nonempty subset of
$\{1,\ldots,d\}$. Again, the order is defined by containment and the
dimension of a face $X$ is $|X|-1$.
Thus, $\beta_d$ has $2^{k+1}{d\choose k+1}$ faces of dimension $k$.

Two infinite series of embeddable skeletons are well-known:
\begin{enumerate}
\item The skeleton of $\alpha_d(\{0\})=\alpha_d(\{d-1\})$ is the
complete graph $K_{d+1}$, which coincides with $J(d+1,1)$.
\item The skeleton of $\beta_d(\{d-1\})=\beta_d(\{0\})^\ast=\gm_d$
coincides with the hypercube graph $H_d$.
\end{enumerate}

The first of these embeddings can be generalized as follows.

\begin{proposition} \label{Johnson graphs}
If $k\in\{0,\ldots,d-1\}$ then the skeleton of $\alpha_d(\{k\})$ coincides
with $J(d+1,k+1)$.
\end{proposition}

\proof We refer to the discussion of the vertices and edges at the end
of Section \ref{WythoffSection}.
According to that discussion, the vertices of $\alpha_d(\{k\})$ are
the $k$-faces, that is, the subsets of $\{1,\ldots,d+1\}$ of
size $k+1$. Furthermore, the only types, leading to edges, are $k-1$ (if
$k>1$) and $k+1$. This means that two vertices are on an edge if and only
if their symmetric difference, as sets, has size two.\qed

Note that the above isomorphism is not quite new, see for example
\cite[page 18]{pak} and \cite[page 8-9]{vallentin}.

\medskip
The above result explains a number of entries in Tables
\ref{WythoffEmbeddable3} and \ref{WythoffEmbeddable4}. We now turn to
the series showing up in line 14 of Table \ref{WythoffEmbeddable3} and
in line 8 of Table \ref{WythoffEmbeddable4}.

\begin{proposition} \label{nearly hypercube}
The skeleton of $\alpha_d(\{0,d-1\})^\ast$ coincides with $H_{d+1}$ with two
antipodal vertices removed. It is an isometric subgraph of $H_{d+1}$.
\end{proposition}

\proof Again we refer to the discussion in Section
\ref{WythoffSection}.  When $V=\{0,d-1\}$, every one-element subset
$U$ of $\Dl=\{0,\ldots,d-1\}$ has the property that $M([U])=U$. This means
that all elements of $\alpha_d(V)$ (that is, all nonempty proper subsets of
$\{1,\ldots,d+1\}$) are vertices of $\alpha_d(V)^\ast$. This also means that
the types corresponding to edges necessarily have size two. If
$U=\{a,b\}\subset\Dl$ and $a<b$ then $M([U])=\{k~|~a\le k\le b\}$. Therefore,
edges of $\alpha_d(V)^\ast$ are flags, whose type is of the form
$\{k,k+1\}$. Thus, we come to the following description of the skeleton
$\Gm$ of $\alpha_d(V)^\ast$: Its vertices are all nonempty proper subsets of
$\{1,\ldots,d+1\}$; two subsets are adjacent when one of them lies in the
other and their sizes differ by one. This matches the well-known definition
of $H_{d+1}$ as the graph on the set of all subsets of
$\{1,\ldots,d+1\}$. In fact, $\Gm$ is $H_{d+1}$ with two antipodal vertices
($\emptyset$ and the entire $\{1,\ldots,d+1\}$) removed. The last claim is
clear.\qed

\medskip
From Proposition \ref{ProbablyFolflore}, we know that 
$\alpha_d(\{0, \dots, d-1\})$ is embeddable into $H_{ {d+1 \choose 2}}$,
which explains line 16 of Table \ref{WythoffEmbeddable3} and line 4 of Table
\ref{WythoffEmbeddable4}.
Also $\beta_d(\{0, \dots, d-1\})$ is embeddable into $H_{d^2}$,
which explains line 6 of Table \ref{WythoffEmbeddable3}
and line 5 of Table \ref{WythoffEmbeddable4}.
There is an easy connection between those two embeddings, that
is the skeleton $\alpha_{d-1}(\{0, \dots, d-2\})$ is a subcomplex of
the complex $\beta_d(\{0, \dots, d-1\})$.

The dual Wythoff polytope in this proposition is, in fact,
the zonotopal Voronoi polytope of the root lattice $\mathsf{A}_d$.
Note that the polytope $\alpha_d(\{0,\ldots,d-1\})$ is known as the {\em
permutahedron}.
It is the zonotopal Voronoi polytope of the dual root lattice $\mathsf{A}_d^\ast$.
Remark also that $\beta_d(\{0,\ldots,d-1\})$ is not the Voronoi polytope of a
lattice, since its number of vertices, $2^dd!$, is greater than $(d+1)!$.

\medskip

\medskip

The following embedding series leaves trace in Tables \ref{WythoffEmbeddable3} and
\ref{WythoffEmbeddable4} in lines 16 and 7, respectively.

\begin{proposition}\label{equivalenceWythoff}
The skeleton of $\beta_d(\{0, \ldots, d-2\})$ is isomorphic to
the skeleton of $\mathsf{D}_d(\{0, \dots, d-1\})$, which is isometrically
embeddable into $H_{d(d-1)}$.
\end{proposition}
\proof We define the following linear functions:
\begin{equation*}
f_i(x)=x_{i+1}-x_{i+2}\mbox{,~for~}0\leq i\leq d-2\mbox{,~and~}f_{d-1}(x)=x_d.
\end{equation*}
Denote by $H_i$ the hyperplane defined by $f_i(x)=0$ and by $s_i$ the
reflection along the hyperplane $H_i$.
Denote by ${\mathcal S}$ the simplex defined by the inequalities $f_i(x)\geq 0$.
The reflections $s_i$ generate the Coxeter group $\mathsf{B}_d$ of fundamental domain
${\mathcal S}$, whose corresponding regular polytope is $\beta_d$.

One way to compute $\beta_d(\{0,\dots, d-2\})$ is to take a vertex 
$v\in {\mathcal S}\cap H_{d-1}$ with $v\notin H_i$ for $i\leq d-2$.
The polytope obtained as convex hull of the
orbit $\mathsf{B}_d(v)$ is then $\beta_d(\{0,\dots,d-2\})$.
Now, the key argument is that ${\mathcal S}\cup s_{d-1}({\mathcal S})$ is also
a simplex ${\mathcal S}'$ defined by the inequalities 
\begin{equation*}
f'_i(x)=f_i(x)\geq 0\mbox{,~for~}0\leq i\leq d-2\mbox{,~and~}f'_{d-1}(x)=x_{d-1}+x_d\geq 0.
\end{equation*}
Those inequalities define hyperplanes which themselves define
orthogonal reflections and so, a finite Coxeter group named $\mathsf{D}_d$.
The point $v$ lies inside of ${\mathcal S}'$ so, the skeleton
of $\beta_d(\{0, \dots d-2\})$ is, actually, the skeleton
of $\mathsf{D}_d(\{0, \dots, d-1\})$.

The isometric embedding follows from Proposition \ref{ProbablyFolflore} and Table \ref{RegularPolytopesAndCoxeterGroups}. \qed

Another way to see the embedding $\beta_d(\{0, \dots, d-2\})$
from $\beta_d(\{0,\dots, d-1\})$ by {\em projection}, that is removing 
dimensions and considering the obtained graph. The idea here is simply to
remove the coordinates corresponding to the planes $x_i=0$.

\medskip
For $d=3$, the polytope $\beta_d(\{0,\ldots,d-2\})$ is the zonotopal Voronoi
polytope of the lattice $\mathsf{A}_3^\ast$.
For $d=4$ and higher, it is a zonotope but it is 
not the Voronoi polytope of a lattice,
since the number of vertices, $2^{d-1}d!$, is greater than $(d+1)!$.

The examples in lines 8 and 9 of Tables \ref{WythoffEmbeddable3} and
\ref{WythoffEmbeddable4}, respectively, suggest that the skeleton of
$\beta_d(\{0,d-1\})$ might be embeddable in a half-cube for all $d$. The next
proposition demonstrates that the actual situation is somewhat more
complicated. We first need to recall some further concepts.

Suppose $\Gm$ is a graph and $\phi$ is a mapping from $\Gm$ to a hypercube
$H_m$. We say that $\phi$ is an {\em embedding with scale $\lm$}
if for all vertices $x,y\in\Gm$ the 
distance in $H_m$ between $\phi(x)$ and $\phi(y)$ (the Hamming distance)
coincides with $\lm d_\Gm(x,y)$. Clearly, isometric embeddings in a
hypercube are scale $1$ embeddings, while isometric embeddings in a
half-cube are scale $2$ embeddings. A graph is an {\em $\ell_1$-graph}
if it has a scale $\lm$ embedding into a hypercube for some $\lm$.
A finite rational-valued 
metric embeds isometrically into some space $\ell_{1}^{k}$ if and only if 
it is scale $\lm$ embeddable into $H_m$ for some $\lm$ and $m$.
See \cite{AD80,Shp,DS1,book1,book2}
for details on $\ell_1$-embedding.

\begin{proposition}
The skeleton of $\beta_d(\{0,d-1\})$ is an $\ell_1$-graph for all
$d$. However, if $d>4$, it is not an isometric subgraph of a half-cube.
\end{proposition}

\proof Let $\Gm$ be the above skeleton graph.  The vertices of this graph
can be identified with all tuples of the form $(k;\eps_1,\ldots,\eps_d)$,
where $1\le k\le d$ and the signs $\eps_i$ are arbitrary. Thus, there are
$d2^d$ vertices. The edges of $\Gm$ arise in two ways: (1) We have that
$(i;\eps_1,\ldots,\eps_d)$ is adjacent to $(j;\eps_1,\ldots,\eps_d)$ for
all $i\ne j$. (2) We also have that $(i;\eps_1,\ldots,\eps_d)$ is adjacent
to $(i;\eps_1,\ldots,\eps_{j-1},-\eps_j,\eps_{j+1}\ldots,\eps_d)$, again
for all $i\ne j$.

Let $\Gm_1$ be the graph whose vertices are all tuples
$(\eps_1,\ldots,\eps_d)$, $\eps_i=\pm 1$, and where two tuples are adjacent
whenever they differ in just one entry. Let $\Gm_2$ be the graph whose
vertices are $\pm k$, $1\le k\le d$, and where vertices $s$ and $t$ are
adjacent whenever $|s|\ne|t|$.  It is clear that $\Gm_1$ is isomorphic to
the hypercube $H_d$, while $\Gm_2$ is isomorphic to the hyperoctahedron
graph $K_{d\times 2}$ (complete multipartite graph with $d$ parts of size
two; also known as the cocktail-party graph). Mapping the vertex
$(k;\eps_1,\ldots,\eps_d)$ of $\Gm$ to the ordered pair
$((\eps_1,\ldots,\eps_d),\eps_kk)$ defines an embedding $\phi$ of $\Gm$
into the Cartesian product graph $\Gm_1\times\Gm_2$. It is easy to see that
this embedding is isometric. Since $\Gm$ projects surjectively onto both
$\Gm_1$ and $\Gm_2$, we can now determine the Graham-Winkler Cartesian product
graph for $\Gm$ (cf. \cite{GW}). Namely, that Cartesian product graph has $d$
complete graphs of size two and the cocktail-party graph $\Gm_2$ as its
factors. (We assume that $d>2$.) It follows from \cite{Shp} that $\Gm$ has
a scale $\lm$ embedding in a hypercube if and only if every factor has.
In our case, every factor is an $\ell_1$-graph, hence $\Gm$ is an
$\ell_1$-graph, too.
Furthermore, the cocktail-party graph $K_{d\times 2}$ with $d>4$
requires $\lm>2$, which proves the second claim of the proposition.\qed

\medskip
The infinite series exhibited in this section explain a majority of the
examples from Tables \ref{WythoffEmbeddable3} and \ref{WythoffEmbeddable4},
including, in fact, {\em all} examples from Table
\ref{WythoffEmbeddable4}. This allows us to make the following conjecture.

\begin{conjecture}\label{conjecture}
If $\Gm$ is the skeleton of the Wythoffian $\cK(V)$ or of the dual
Wythoffian $\cK(V)^\ast$, where $\cK$ is a regular convex polytope, and
$\Gm$ is isometrically embeddable in a half-cube, then $\Gm$ can be found
either in Table \ref{WythoffEmbeddable3}, Table \ref{WythoffEmbeddable4},
or in one of the infinite series discussed in this section.
\end{conjecture}

\section{From the hypercubes to the cubic lattices?}

The above conjecture shows one direction of possible further
research. Another possibility is extending the results of this paper to
cover the case of infinite regular polytopes, that is, regular partitions
of the Euclidean and hyperbolic space. In this section we briefly discuss
what is known about the easier Euclidean case.

In the infinite case, instead of embedding the skeleton graphs up to scale
into hypercubes $H_m$, we embed them into the $m$-dimensional cubic lattice
$Z_m$ (including $m=\infty$) taken with its metric $\ell_1$.
Notice that this is a true generalization, because, according to \cite{AD80},
a finite metric that can be embedded into a cubic lattice, can also be
embedded into a hypercube.

All regular partitions of Euclidean $d$-space ($d$ finite) are known
\cite{Cox73}. They consist of one infinite series $\delta_d=\delta_d^\ast$, which
is the partition into the regular $d$-dimensional cubes, two
$2$-dimensional ones, $(3^6)$ (partition into regular triangles) and
$(6^3)=(3^6)^\ast$ (partition into regular $6$-gons), and two
$4$-dimensional ones, $h\delta_4$ (partition into $4$-dimensional
hyperoctahedra) and $h\delta_4^\ast$ (partition $24$-cells). Notice that the
latter two partitions are the Delaunay and Voronoi partitions associated
with the lattice $\mathsf{D}_4$. In particular, below we use the notation $Vo(\mathsf{D}_4)$
in place of $h\delta_4^\ast$.

In the following table we give a complete list of Wythoffians of regular
partitions of the Euclidean plane. We use the classical notation for the
vertex-transitive partition of the Euclidean plane; namely, each partition
is identified by its {\em type}, listing clock-wise the sizes of the
faces containing a fixed vertex. In particular, the regular partitions of
the Euclidean plane are $(4^4)=\delta_2$, $(3^6)$, and $(6^3)=(3^6)^\ast$. In
the second column we indicate the embedding. We put $Z_m$ for an embedding
with scale one and $\frac{1}{2}Z_m$ for an embedding with scale two.
\begin{table}[ht]
\begin{center}
\begin{tabular}{||c|c||}
\hline
\hline
Wythoffian                          & embedding\\
\hline
$\delta_2=\delta_2(\{0\})=\delta_2(\{1\})=\delta_2(\{2\})=\delta_2(\{0,2\})$ & $Z_2$\\
$(3^6)=(3^6)(\{0\})$                & $\frac{1}{2}Z_3$\\
$(6^3)=(3^6)(\{2\})=(3^6)(\{0,1\})$ & $Z_3$\\
\hline
$(4.8^2)=\delta_2(\{0,1\})=\delta_2(\{1,2\}=\delta_2(\{0,1,2\}$ & $Z_4$\\
$(4.6.12)=(3^6)(\{0,1,2\})$         & $Z_6$\\
$(3.4.6.4)=(3^6)(\{0,2\})$          & $\frac{1}{2}Z_3$\\
$(3.6.3.6)^\ast=(3^6)(\{1\})^\ast$  & $Z_3$\\
$(3.12^2)^\ast=(3^6)(\{1,2\})^\ast$ & $\frac{1}{2}Z_{\infty}$ \\
\hline
\hline
\end{tabular}
\end{center}
\caption{Embeddable Wythoffian cases for plane partitions.}
\label{WythoffEmbeddable3i}
\end{table}

All Archimedean Wythoffians or their dual, which are not mentioned in Table
\ref{WythoffEmbeddable3i}, are non-embeddable and, moreover, they do not
satisfy the $5$-gonal inequality.

In this table we separated the three regular plane partitions from the {\em
Archimedean} ({\it i.e.}, vertex- but not face-transitive) ones. Notice
that, out of the eight Archimedean partitions, five are
Wythoffians. Missing are partitions $(3^2.4.3.4)$, $(3^3.4^2)$ and
$(3^4.6)$. It turns out that for all regular and Archimedean plane
partitions (and in particular, for all our Wythoffians) exactly one out of
itself and its dual is embeddable. In this respect the situation here
repeats the situation for the Archimedean polyhedra for $d=3$, see Section
\ref{dimension 2 and 3} and Tables 9.1 and 4.1--4.2 in \cite{book2}.

We now turn to the next dimension, $d=3$. Here we identify the Wythoffians
as particular partitions of the Euclidean $3$-space in two ways. First, in
column $2$ we give the number of that partition in the list of 28 regular
and Archimedean partitions of the $3$-space from \cite{book2}. Secondly, we
identify in column 3 the tiles of the partition.  Here, as before, $\beta_3$
and $\gm_3$ are the Octahedron and the Cube, respectively. Also, ``Cbt''
stands for the Cuboctahedron and ``Rcbt'' stands for the
Rhombicuboctahedron. Clearly, ``tr'' stands for ``truncated'' and Prism$_8$
is the regular faced $8$-gonal prism.
In some cases we also indicate the chemical
names of the corresponding partitions. In column 4 we give the details of
the embedding. If the particular Wythoffian is non-embeddable, we put ``non
5-gonal'' in that column to indicate that it fails the $5$-gonal
inequality. The information in column 4 is taken from Table 10.1 from
\cite{book2}.
\begin{table}[ht]
\begin{center}
\begin{tabular}{||c|c|c|c||}
\hline 
\hline 
Wythoffian                                       & no & tiles & embedding\\ 
\hline
$\delta_3=\delta_3(\{0\})=\delta_3(\{3\})=\delta_3(\{0,3\})$ & 1 & $\gm_3$ & $Z_3$\\
$\delta_3(\{1,2\})=Vo(\mathsf{A}_3^\ast)$                    & 2 & tr $\beta_3$ & $Z_6$\\
$\delta_3(\{0,1,2\})=\delta_3(\{1,2,3\})$=zeolit Linde A & 16 & $\gm_3$, tr $\beta_3$, tr Cbt & $Z_9$\\
$\delta_3(\{0,1,2,3\})$=zeolit $\rho$               & 9 & Prism$_8$, tr Cbt & $Z_9$\\
\hline
$\delta_3(\{1\})=\delta_3(\{2\})=De(J-complex)$        & 8 & $\beta_3$, Cbt & non 5-gonal\\ 
$\delta_3(\{0,1\})=\delta_3(\{2,3\})$=boride $CaB_6$   & 7 & $\beta_3$, tr $\gm_3$ & non 5-gonal\\
$\delta_3(\{0,2\})=\delta_3(\{1,3\})$                  & 18 & $\gm_3$, Cbt, Rcbt & non 5-gonal\\
$\delta_3(\{0,1,3\})=\delta_3(\{0,2,3\})$=selenide $Pd_{17}Se{15}$ & 23 & $\gm_3$, Prism$_8$, tr $\gm_3$, Rbct & non 5-gonal\\
\hline
\hline
\end{tabular}
\end{center}
\caption{Wythoffians of regular partitions of the $3$-space.}
\label{WythoffEmbeddable4i}
\end{table}

As Table \ref{WythoffEmbeddable4i} indicates, only eight out of 28 regular
and Archimedean partitions of the $3$-space arise as the Wythoffians of the
cubic partition $\delta_3$.

Finally, in Table \ref{WythoffEmbeddable5i} we collected some information
about the dimensions $d\ge 4$.

\begin{table}[ht]
\begin{center}
\begin{tabular}{||c|c|c||}
\hline
\hline
Wythoffian                                       & tiles & embedding\\
\hline
$\delta_d=\delta_d(\{0\})=\delta_d(\{d\})=\delta_d(\{0,d\})$ & $\gm_d$ & $Z_d$\\
$\delta_d(\{0,1\})$=tr $\delta_d$                      & $\beta_d$, tr $\gm_d$ & non 5-gonal\\
\hline
$Vo(\mathsf{D}_4)=Vo(\mathsf{D}_4)(\{0\}$                          & $24-cell$ & non 5-gonal\\
$Vo(\mathsf{D}_4)^\ast=Vo(\mathsf{D}_4)(\{4\})$                    & $\beta_4$ & non 5-gonal\\
$Vo(\mathsf{D}_4)(\{1\}=Med(Vo(\mathsf{D}_4))$                     & $\gm_4$, $Med(24-cell)$ & non 5-gonal\\
$Vo(\mathsf{D}_4)(\{0,1\}$=tr $Vo(\mathsf{D}_4)$                   & $\gm_4$, tr $24-cell$ & $Z_{12}$\\
\hline
\hline
\end{tabular}
\end{center}
\caption{Some Wythoffians of regular partitions of the $d$-space, $d\ge 4$.}
\label{WythoffEmbeddable5i}
\end{table}

Notice that again, as in Table \ref{WythoffEmbeddable4}, few Wythoffians
for $d=3$ possess embeddings. This gives hope that there is only a small
number of infinite series of embeddings in the Euclidean case.  One of the
infinite series is shown in line 1 of Table
\ref{WythoffEmbeddable5i}.

\begin{proposition}\label{affinecase}
(i) The skeleton graph of the Wythoffian flag complex $\delta_d(\{0,\ldots,d\})$
is isometrically embeddable in $Z_{d^2}$.

(ii) The skeleton graph of the Wythoffian
$\delta_d(\{0,\ldots,d-1\})=\delta_d(\{1,\ldots,d\})$ is also isometrically
embeddable in $Z_{d^2}$.
\end{proposition}
\proof The situation is completely analogous to the one for
$\beta_d(\{0,\dots, d-1\})$ and $\beta_d(\{0, \dots, d-2\})$
of Proposition \ref{equivalenceWythoff}.
We simply add the function $f_d(x)=1-x_1$ and the Coxeter
groups $\mathsf{B}_d$ and $\mathsf{D}_d$ are replaced,
respectively, by $\widetilde{\mathsf{C}_d}$
and $\widetilde{\mathsf{B}_d}$ (see \cite[page~96]{humphreyscoxeter}).
The number of classes of parallel hyperplanes in $\widetilde{G}$ is equal
to the number of hyperplanes in $G$. The result follows from Table
\ref{RegularPolytopesAndCoxeterGroups}. \qed

We think that those tilings are zonotopal but we did not check it.

It appears (see line 4 of Table \ref{WythoffEmbeddable3i} and line 2 of
Table \ref{WythoffEmbeddable4i}) that $\delta_d(\{1,2\})$ may have an
embedding for all $d$. In this case, however, we are reluctant to formulate
an exact conjecture. Perhaps, the situation will be more clear when the
case $d=4$ is completed.

We have already pointed out that only eight out of 28 regular and
Archimedean partitions of the Euclidean $3$-space are Wythoffians of
$\delta_3$. This indicates that, maybe, one needs to derive Wythoffians from a
larger class of Euclidean partitions. The obvious candidates are the
Delaunay and Voronoi partitions of interesting Euclidean lattices, in
particular, the root lattices. For the case of such lattices themselves
(that is, for $V=\{0\}$ or $\{d\}$), see Chapter 11 of \cite{book2}.
The zonotopal embeddings of $Vo(\mathsf{A}_d)$ in $Z_{d+1}$
and $Vo(\mathsf{A}_d^\ast)$ in $Z_{{d+1\choose 2}}$,
correspond to the zonotopal embeddings of the corresponding tiles
from Proposition \ref{ProbablyFolflore}.

It is, of course, also very interesting to consider the Wythoffians of the
regular partitions of the {\em hyperbolic} $d$-space.  In fact, in
\cite{DS} (see also Chapter 3 of \cite{book2}), the embeddability was
decided for any regular tiling $P$ of the $d$-sphere, Euclidean $d$-space,
hyperbolic $d$-space or Coxeter's regular hyperbolic honeycomb (with
infinite or star-shaped cells or vertex figures). The large program will be
to generalize it for all Wythoffians of such general $P$.

\end{document}